\documentclass[12pt,oneside]{amsart}

\usepackage[a4paper,body={14.6cm, 23.5cm}, mag=1000]{geometry}
\usepackage{amssymb}
\usepackage{amsthm}

\theoremstyle{plain}
\newtheorem{theorem}{Theorem}[section]

\newtheorem{prop}{Proposition}[section]
\newtheorem{corollary}{Corollary}[section]

\begin{document}

\title[TWO-SIDED BOUNDS FOR VOLUMES]
{Two-sided bounds for the volume of right-angled hyperbolic polyhedra}

\author{Du\v{s}an Repov\v{s}}
\address{Faculty of Mathematics and Physics, and Faculty
of Education, University of Ljubljana, P.O.Box 2964,
SI-1001 Ljubljana, Slovenia}
\email{dusan.repovs@guest.arnes.si}

\author{Andrei Vesnin}
\address{Sobolev Institute of Mathematics, pr. ak. Koptyuga~4,
Novosibirsk, 630090, Russia \and Department of Mathematics, Omsk
State Technical University, pr. Mira~11, Omsk, 644050, Russia}
\email{vesnin@math.nsc.ru}


\thanks{The first author was supported in part by the Slovenian
Research Agency grants P1-0292-0101, J1-9643-0101 and J1-2057-0101.
The second author was supported in part by the Russian Foundation
for Basic Research grant 09-01-00255 and by the SO RAN -- UrO RAN
grant.}


\begin{abstract} For a compact right-angled
polyhedron $R$ in $\mathbb H^3$ denote by $\operatorname{vol} (R)$
the volume and by $\operatorname{vert} (R)$ the number of vertices.
Upper and lower bounds for $\operatorname{vol} (R)$ in terms of
$\operatorname{vert} (R)$ were obtained in \cite{A09}. Constructing
a 2-parameter family of polyhedra, we show that the asymptotic upper
bound $5 v_3 / 8$, where $v_3$ is the volume of the ideal regular
tetrahedron in $\mathbb H^3$, is a double limit point for ratios
$\operatorname{vol} (R) / \operatorname{vert} (R)$. Moreover, we
improve the lower bound in the case $\operatorname{vert} (R)
\leqslant 56$.
\end{abstract}

\date{\today}

\maketitle

\section{Right-angled polyhedra in $\mathbb H^3$.}

In any space, right-angled polyhedra are very convenient to serve as
``building blocks'' for various geometric constructions. In
particular, they have several interesting properties in hyperbolic
3-space $\mathbb H^3$. One can try to obtain a hyperbolic 3-manifold
using a right-angled polyhedron as its fundamental polyhedron. Or,
one can construct a hyperbolic 3-manifold in such a way that its
fundamental group is a torsion-free subgroup of the Coxeter group,
generated by reflections across the faces of a right-angled
polyhedron \cite{V87}. Below we consider only compact polyhedra,
which do not admit ideal vertices.

We start by recalling two nice recent results. Inoue \cite{I08}
introduced two operations on right-angled polyhedra called
\emph{decomposition} and \emph{edge surgery}, and proved that
L\"obell polyhedra (which will be a subject of discussion below) are
universal in the following sense:

\begin{theorem} \cite[Theorem~9.1]{I08} \label{theorem1.1}
Let $P_0$ be a right-angled hyperbolic polyhedron. Then there exists
a sequence of disjoint unions of right-angled hyperbolic polyhedra
$P_1, \ldots, P_k$ such that for $i=1, \ldots, k$, $P_i$ is obtained
from $P_{i-1}$ by either a decomposition or an edge surgery, and
$P_k$ is a set of L\"obell polyhedra. Furthermore,
$$
\operatorname{vol} (P_0) \geqslant \operatorname{vol} (P_1)
\geqslant \operatorname{vol} (P_2) \geqslant \ldots \geqslant
\operatorname{vol} (P_k) .
$$
\end{theorem}

Atkinson \cite{A09} estimated the volume of a right-angled
polyhedron in terms of the number of its vertices as follows:

\begin{theorem} \cite[Theorem~2.3]{A09} \label{theorem1.2}
If $P$ is a compact right-angled hyperbolic polyhedron with $V$
vertices, then
$$
(V-2) \cdot \frac{v_8}{32} \leqslant \operatorname{vol} (P) < (V-10)
\cdot \frac{5 v_3}{8},
$$
where $v_8$ is the volume of a regular ideal octahedron, and $v_3$
is the volume of a regular ideal tetrahedron. There is a sequence of
compact polyhedra $P_i$, with $V_i$ vertices such that
$\operatorname{vol} (P_i) / V_i$ approaches $5 v_3 / 8$ as $i$ goes
to infinity.
\end{theorem}

A family of polyhedra $P_i$ suggested by Atkinson is described  in
the proof of \cite[Prop.~6.4]{A09}.

In this note we will demonstrate that L\"obell polyhedra can serve
as a suitable family realizing the upper bound. Thus these polyhedra
play an important role not only in Theorem~\ref{theorem1.1}, but
also in Theorem~\ref{theorem1.2}.

Let us denote by $\operatorname{vert} (R)$ the number of vertices of
a right-angled polyhedron $R$. In this note we prove that $5 v_3 /
8$ is a double limit point in the sense that it is the limit point
of limit points for ratios $\operatorname{vol} (R) /
\operatorname{vert} (R)$.

\begin{theorem} \label{theorem1.3}
For any integer $k \geqslant 1$ there exists a series of compact
right-angled polyhedra $R_k(n)$ in $\mathbb H^3$ such that
$$
\lim_{n \to \infty} \frac{\operatorname{vol} (R_k(n))}
{\operatorname{vert} (R_k(n))} = \frac{k}{k+1} \cdot \frac{5 v_3}{8}
.
$$
\end{theorem}

As one will see from the proof, $R_1(n)$ are L\"obell polyhedra and
$R_k(n)$ for $k>1$ are towers of them.

Moreover, in Corollary~\ref{corollary4.3} we improve the lower
estimate from Theorem~\ref{theorem1.2} in the case
$\operatorname{vert} (R) \leqslant 56$.

\section{L\"obell polyhedra and manifolds.}

We introduced L\"obell polyhedra in \cite{V87} as a generalization
of a right-angled $14$-hedron used in \cite{L31}.

\begin{figure}[ht]
\begin{center}
\setlength{\unitlength}{0.36mm}
\begin{picture}(240,100)(-120,10)
\thicklines \put(-120,0){\begin{picture}(120,120)
\put(40,50){\line(0,1){20}} \put(40,70){\line(2,1){20}}
\put(60,80){\line(2,-1){20}} \put(80,70){\line(0,-1){20}}
\put(80,50){\line(-2,-1){20}} \put(60,40){\line(-2,1){20}}
\put(40,50){\line(-1,-1){10}} \put(40,70){\line(-1,1){10}}
\put(60,80){\line(0,1){10}} \put(80,70){\line(1,1){10}}
\put(80,50){\line(1,-1){10}} \put(60,40){\line(0,-1){10}}
\put(40,25){\line(-2,3){10}} \put(40,25){\line(4,1){20}}
\put(20,60){\line(1,2){10}} \put(20,60){\line(1,-2){10}}
\put(40,95){\line(4,-1){20}} \put(40,95){\line(-2,-3){10}}
\put(80,95){\line(-4,-1){20}} \put(80,95){\line(2,-3){10}}
\put(100,60){\line(-1,2){10}} \put(100,60){\line(-1,-2){10}}
\put(80,25){\line(-4,1){20}} \put(80,25){\line(2,3){10}}
\put(30,10){\line(-3,5){30}} \put(0,60){\line(3,5){30}}
\put(30,110){\line(1,0){60}} \put(90,110){\line(3,-5){30}}
\put(120,60){\line(-3,-5){30}} \put(90,10){\line(-1,0){60}}
\put(30,10){\line(2,3){10}} \put(0,60){\line(1,0){20}}
\put(30,110){\line(2,-3){10}} \put(90,110){\line(-2,-3){10}}
\put(120,60){\line(-1,0){20}} \put(90,10){\line(-2,3){10}}
\end{picture} }
\put(20,0){\begin{picture}(100,120) \put(50,40){\line(-2,1){20}}
\put(30,50){\line(1,4){5}} \put(35,70){\line(1,0){30}}
\put(65,70){\line(1,-4){5}} \put(70,50){\line(-2,-1){20}}
\put(50,40){\line(0,-1){10}} \put(30,50){\line(-2,-1){10}}
\put(35,70){\line(-1,2){5}} \put(65,70){\line(1,2){5}}
\put(70,50){\line(2,-1){10}} \put(30,25){\line(-1,2){10}}
\put(30,25){\line(4,1){20}} \put(30,25){\line(-2,-3){10}}
\put(15,65){\line(1,1){15}} \put(15,65){\line(1,-4){5}}
\put(15,65){\line(-3,1){15}} \put(50,90){\line(2,-1){20}}
\put(50,90){\line(-2,-1){20}} \put(50,90){\line(0,1){20}}
\put(85,65){\line(-1,-4){5}} \put(85,65){\line(-1,1){15}}
\put(85,65){\line(3,1){15}} \put(70,25){\line(-4,1){20}}
\put(70,25){\line(1,2){10}} \put(70,25){\line(2,-3){10}}
\put(20,10){\line(2,3){10}} \put(0,70){\line(3,-1){15}}
\put(50,110){\line(0,-1){20}} \put(100,70){\line(-3,-1){15}}
\put(80,10){\line(-2,3){10}} \put(20,10){\line(-1,3){20}}
\put(0,70){\line(5,4){50}} \put(50,110){\line(5,-4){50}}
\put(100,70){\line(-1,-3){20}} \put(80,10){\line(-1,0){60}}
\end{picture} }
\end{picture}
\end{center}
\caption{Polyhedra $R(6)$ and $R(5)$.}
\label{fig1}
\end{figure}
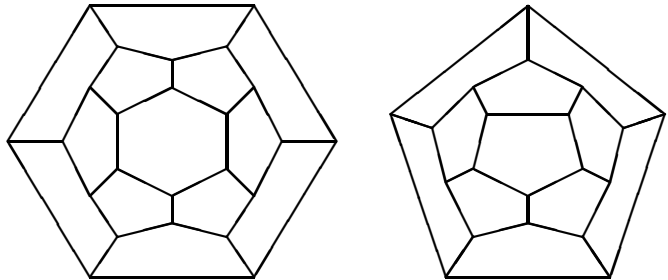

Recall that in order to give a positive answer to the question of
the existence of ``Clifford-Klein space forms'' (that is, closed
manifolds) of constant negative curvature, L\"obell~\cite{L31}
constructed in 1931 the first example of a closed orientable
hyperbolic 3-manifold. This manifold was obtained by gluing together
eight copies of the right-angled $14$-faced polytope (denoted below
by $R(6)$ and shown in Fig.~\ref{fig1}) with an upper and a lower
basis both being regular hexagons, and a lateral surface given by
$12$ pentagons, arranged similarly as in the dodecahedron.
Obviously, $R(6)$ can be considered as a generalization of a
right-angled dodecahedron in the way of replacing basis pentagons to
hexagons.

As shown in \cite{V87}, the dodecahedron and $R(6)$ are part of a
larger family of polyhedra. For each $n\geqslant 5$ we consider the
right-angled polyhedron $R(n)$ in $\mathbb H^3$ with $(2n+2)$ faces,
two of which (viewed as the upper and  lower bases) are regular
$n$-gons, while the lateral surface is given by $2n$ pentagons,
arranged as one can easily imagine. Note that $R(5)$ is the
right-angled dodecahedron (see Fig.~\ref{fig1}). Existence of
polyhedra $R(n)$ in $\mathbb H^3$ can be easy checked by involving
Andreev's theorem \cite{An}.

An algebraic approach suggested in \cite{V87} admits a construction
of both orientable and non-orientable closed hyperbolic
$3$-manifolds from eight copies of any bounded right-angled
hyperbolic polyhedron. More exactly, any coloring of the faces of a
right-angled polyhedron by four colors so that no two faces of the
same color share an edge encodes a torsion-free subgroup of
orientation preserving isometries which is a subgroup of the
polyhedral Coxeter group of index eight. Thus, any four-coloring
encodes an orientable hyperbolic 3-manifold obtained from eight
copies of a right-angled polyhedron. This approach also allows one
to construct non-orientable hyperbolic 3-manifolds, but in this case
five to seven colors are needed.

It was mentioned in \cite{V87} that the manifold constructed by
L\"obell can be encoded by some four-coloring of $R(6)$, and it was
shown how to construct concrete orientable and non-orientable
manifolds using eight copies of $R(n)$ for any $n\geqslant 5$.
Closed orientable hyperbolic 3-manifolds encoded by four-colorings
of $R(n)$, $n\geqslant 5$, were called \emph{L\"obell manifolds}.
(Observe that for each $n$ number of such manifolds do not need to
be unique.) Polyhedra $R(n)$ can be naturally referred as
\emph{L\"obell polyhedra}.

Various properties of L\"obell manifolds were intensively studied:
the volume formulae were obtained in \cite{MedV07} and \cite{V98},
invariant trace fields for fundamental groups and their
arithmeticity were numerically calculated in \cite{Roeder}, many of
L\"obell manifolds were obtained in \cite{MedV99} as two-fold
branched coverings of the 3-sphere, and two-sided bounds for
complexity of L\"obell manifolds were done in \cite{MPV09}.

Since Lobachevsky's 1832 paper, the following \emph{Lobachevsky
function} has traditionally been used in volume formulae for
hyperbolic polyhedra
$$
\Lambda(x) = - \int\limits_0^x \log | 2 \sin (t) | \, {\rm d} t.
$$

The volume formula for L\"obell manifolds established in \cite{V98}
implies the following formula for $\operatorname{vol} R(n)$, since
any L\"obell manifolds indexed by $n$ is glued by isometries from
eight copies of $R(n)$:

\begin{theorem} \label{theorem2.1}
For all $n\geqslant 5$ we have
$$
\operatorname{vol} (R(n))
=  \frac{n}{2} \left( 2 \Lambda (\theta_n)  +  \Lambda \left(
\theta_n + \frac{\pi}{n} \right)  +  \Lambda \left( \theta_n -
\frac{\pi}{n} \right)  +  \Lambda \left( \frac{\pi}{2} - 2 \theta_n
\right) \right) ,
$$
where $$\theta_n  =  \frac{\pi}{2}  -  \arccos  \left(\frac{1}{2
\cos ( \pi / n )} \right).$$
\end{theorem}

It is easy to check that $\theta_n \to \pi/6$ and
$\frac{\operatorname{vol} R(n)}{n} \to \frac{5 v_3}{4}$ as $n \to
\infty$. Here we use that $v_3 = 3 \Lambda(\pi/3) = 2
\Lambda(\pi/6)$. Moreover, the asymptotic behavior of volumes of
L\"obell manifolds was established in \cite[Prop.~2.10]{MPV09}. This
implies trivially the description of the asymptotic behavior of
$\operatorname{vol} (R(n))$ as $n$ tends to infinity.

\begin{prop} \label{prop2.1}
The following inequalities hold for sufficiently large~$n$:
$$
\frac{5v_3}{4} \cdot n - \frac{17 v_3}{2n}  <
\operatorname{vol} (R(n))
< \frac{5v_3}{4} \cdot n.
$$
\end{prop}

Since $\operatorname{vert} (R(n)) = 4 n$, we get

\begin{corollary} \label{corollary2.1}
The following inequalities hold for sufficiently large~$n$:
$$
\frac{5v_3}{16} - \frac{17 v_3}{8n^2}  <
\frac{\operatorname{vol} (R(n))}{\operatorname{vert} (R(n))}
< \frac{5 v_3}{16} .
$$
\end{corollary}

\section{Proof of Theorem~\ref{theorem1.3}.}

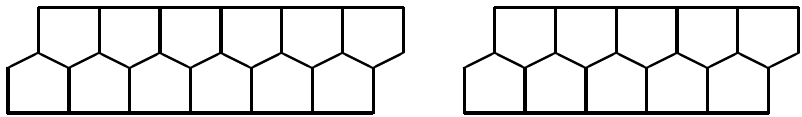
\begin{figure}[ht]
\begin{center}
\setlength{\unitlength}{0.4mm}
\begin{picture}(240,35)(-120,25)
\thicklines \put(-130,0){\begin{picture}(130,80)
\put(0,25){\line(1,0){120}} \put(10,60){\line(1,0){120}}
\multiput(0,25)(20,0){7}{\line(0,1){15}}
\multiput(10,60)(20,0){7}{\line(0,-1){15}}
\multiput(0,40)(20,0){7}{\line(2,1){10}}
\multiput(10,45)(20,0){6}{\line(2,-1){10}}
\end{picture} }
\put(20,0){\begin{picture}(110,80) \put(0,25){\line(1,0){100}}
\put(10,60){\line(1,0){100}}
\multiput(0,25)(20,0){6}{\line(0,1){15}}
\multiput(10,60)(20,0){6}{\line(0,-1){15}}
\multiput(0,40)(20,0){6}{\line(2,1){10}}
\multiput(10,45)(20,0){5}{\line(2,-1){10}}
\end{picture} }
\end{picture}
\end{center}
\caption{Polyhedra $R(6)$ and $R(5)$.} \label{fig2}
\end{figure}

We will use L\"obell polyhedra $R(n)$ as building blocks to
construct right-angled polyhedra with necessary properties. Let us
present polyhedra $R(n)$ by their lateral surfaces as it is done in
Fig.~\ref{fig2} for polyhedra $R(6)$ and $R(5)$, keeping in mind
that left and right sides are glued together.

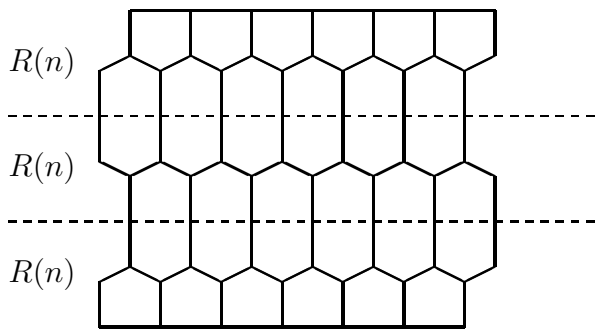
\begin{figure}[ht]
\begin{center}
\setlength{\unitlength}{0.4mm}
\begin{picture}(130,105)(0,0)
\thicklines
\put(0,0){\line(1,0){120}}
\multiput(0,0)(20,0){7}{\line(0,1){15}}
\multiput(10,35)(20,0){7}{\line(0,-1){15}}
\multiput(0,15)(20,0){7}{\line(2,1){10}}
\multiput(10,20)(20,0){6}{\line(2,-1){10}}
{\thinlines \multiput(-30,35)(5,0){40}{\line(1,0){2.5}}}
\put(-30,15){$R(n)$}
\multiput(10,35)(20,0){7}{\line(0,1){15}}
\multiput(0,70)(20,0){7}{\line(0,-1){15}}
\multiput(10,50)(20,0){6}{\line(2,1){10}}
\multiput(0,55)(20,0){7}{\line(2,-1){10}}
{\thinlines \multiput(-30,70)(5,0){40}{\line(1,0){2.5}}}
\put(-30,50){$R(n)$}
\put(-30,85){$R(n)$}
\multiput(0,70)(20,0){7}{\line(0,1){15}}
\multiput(10,105)(20,0){7}{\line(0,-1){15}}
\multiput(0,85)(20,0){7}{\line(2,1){10}}
\multiput(10,90)(20,0){6}{\line(2,-1){10}}
\put(10,105){\line(1,0){120}}
\end{picture}
\end{center}
\caption{Polyhedron $R_3(6)$.} \label{fig3}
\end{figure}

For integer $k \geq 1$ denote by $R_k (n)$ the polyhedron
constructed from $k$ copies of $R(n)$ gluing them along $n$-gonal
faces similar to a tower. In particular, $R_1 (n) = R(n)$. The
polyhedron $R_3(6)$ is presented in Fig.~\ref{fig3}.

Obviously, $R_k(n)$ is a right-angled polyhedron with $n$-gonal top
and bottom and the lateral surface formed by $2n$ pentagons and
$(k-1)n$ hexagons.

Since $\operatorname{vol} (R_k(n)) = k \cdot \operatorname{vol}
(R(n))$, Proposition~\ref{prop2.1} implies that for sufficiently
large $n$
$$
k  \cdot \frac{5v_3}{4} \cdot n - k \cdot \frac{17 v_3}{2 n}  <
\operatorname{vol} (R_k(n)) < k \cdot \frac{5v_3}{4} \cdot n.
$$
Since $\operatorname{vert} R_k(n) = (2k+2) n$, we obtain
$$
\frac{k}{k+1} \cdot \frac{5v_3}{8}  - \frac{k}{k+1} \cdot \frac{17
v_3}{4 n^2} < \frac{\operatorname{vol} (R_k(n))}{\operatorname{vert}
(R_k(n))} < \frac{k}{k+1} \cdot \frac{5 v_3}{8} .
$$

Thus family of right-angled polyhedra $R_k(n)$ is such that for any
integer $k \geqslant 1$
$$
\lim_{n \to \infty} \frac{\operatorname{vol}
(R_k(n))}{\operatorname{vert} (R_k(n))} = \frac{k}{k+1} \cdot
\frac{5 v_3}{8},
$$
and the upper bound $5 v_3 / 8$ is a double limit point in the sense
that it is the limit of above limit points as $k \to \infty$:
$$
\lim_{k, n \to \infty} \frac{\operatorname{vol}
(R_k(n))}{\operatorname{vert} (R_k(n))} = \frac{5 v_3}{8}.
$$
Thus, the theorem is proved. \qed

\section{Other volume estimates.}

Since 1-skeleton of a right-angled compact hyperbolic polyhedron
$P$ is a trivalent plane graph, one can easy see that
Euler formula for a polyhedron implies
$$
V = 2F - 4 ,
$$
where $V$ is number of vertices of $P$ and $F$ is number of
its faces. Moreover, Euler formula implies also that
$P$ has at least $12$ faces (this smallest number of faces
corresponds to a dodecahedron).
Thus, Theorem~\ref{theorem1.2} implies the
following result.

\begin{corollary} \label{corollary4.1}
If $P$ is a compact right-angled hyperbolic polyhedron with $F$
faces, then
$$
(F-3) \cdot \frac{v_8}{16} \leqslant \operatorname{vol} (P) < (F -
7) \cdot \frac{5 v_3}{4} .
$$
\end{corollary}

We recall that constants $v_3$ and $v_8$ are
$$
v_3 = 3 \, \Lambda ( \pi/3 )  = 1.0149416064096535\ldots
$$
and
$$
v_8 = 8 \, \Lambda ( \pi / 4 ) = 3.663862376708876\dots .
$$

Since a right-angled hyperbolic $n$-gon has area $\pi / 2 \cdot ( n
- 4)$, the lateral surface area of a compact hyperbolic right-angled
polyhedron $P$ with $F$ faces is equal to $\pi \cdot (F - 6)$. Thus,
Corollary~\ref{corollary4.1} implies the following result.

\begin{corollary} \label{corollary4.2}
If $P$ is a compact right-angled hyperbolic polyhedron with lateral
surface area $S$, then
$$
(S / \pi + 3 ) \cdot \frac{v_8}{16} \leqslant \operatorname{vol} (P)
< (S / \pi - 1 ) \cdot \frac{5 v_3}{4} .
$$
\end{corollary}

Observe, that Theorem~\ref{theorem2.1} can be used to show that the
volume function $\operatorname{vol} R(n)$ is a monotonic increasing
function of $n$ (see \cite{I08} and \cite{MPV09} for proofs), and to
calculate volumes of L\"obell polyhedra. In particular,
$$
\operatorname{vol} R(5) = 4.306\ldots, \qquad
\operatorname{vol} R(6) = 6.023\ldots, \qquad
\operatorname{vol} R(7) = 7.563\ldots.
$$
Together with Theorem~\ref{theorem1.1} it gives that the
right-angled hyperbolic polyhedron of smallest volume is $R(5)$ (a
dodecahedron) and the second smallest is $R(6)$. Thus, if a compact
right-angled hyperbolic polyhedron $P$ is differ of a dodecahedron,
then
$$
\operatorname{vol} (P) \geqslant 6.023\ldots .
$$
Thus, we get the following

\begin{corollary} \label{corollary4.3}
If $P$ is a compact right-angled hyperbolic polyhedron different
than a dodecahedron, having $V$ vertices and $F$ faces. Then
$$
\operatorname{vol} (P) \geqslant \max \{ (V-2) \cdot \frac{v_8}{32},
\,  6.023\ldots \}
$$
and
$$
\operatorname{vol} (P) \geqslant \max \{ (F-3) \cdot \frac{v_8}{16},
\, 6.023\ldots \} .
$$
\end{corollary}

The estimates from Corollary~\ref{corollary4.3} improve the lower
estimate from Theorem~\ref{theorem1.2} for $V \leqslant 54$ and the
lower estimate from Corollary~\ref{corollary4.1} for $F \leqslant
29$.

\end{document}